\documentstyle[leqno,12pt]{article}
\input{amssym.def}
\input{amssym.tex}

\def\rz{{\rm I \! R}}
\def\cz{{\unitlength1pt\begin{picture}(9,8)
         \put(0,0){{\rm C}}
         \put(3,0.5){\line(0,1){7.5}}
         \end{picture}}}

\def\scz{{\,\rm I \!\!\! C}}
\def\nz{{\Bbb N}}
\def\gz{{\Bbb Z}}
\def\pr{{\rm I \! P}}

\newcommand{\skp}{\hspace{1pt}}
\newcommand{\eps}{\varepsilon}
\newcommand{\ph}{\varphi}
\newcommand{\tr}{{\rm tr\skp}}

\newcommand{\rk}{{\mbox{\rm rk}}}

\newcommand{\Hom}{{\mbox{\rm Hom}\skp}}

\newcommand{\di}{\partial}
\newcommand{\dbar}{{\bar\di}}

\newcommand{\cale}{{\cal E}}
\newcommand{\calf}{{\cal F}}

\newcommand{\cali}{{\cal I}}
\newcommand{\calj}{{\cal J}}
\newcommand{\calm}{{\cal M}}
\newcommand{\calo}{{\cal O}}

\newtheorem{prop}{Proposition}[section]
\newtheorem{theorem}[prop]{Theorem}
\newtheorem{lemma}[prop]{Lemma}

\newtheorem{rem}[prop]{Remark}

\newcommand{\pf}{{\em Proof. }}
\newcommand{\qed}{{\hfill$\Diamond$}\vspace{1.5ex}}

\newcommand{\grass}{{\rm G}(k,n)}
\newcommand{\qrM}{H^*_{[\omega]}(M)}
\newcommand{\qrG}{H^*_{[\omega]}(\grass)}
\newcommand{\qwedge}{\wedge_{\mbox{\scriptsize Q}}}
\newcommand{\bflm}{{\underline{\lambda}}}
\newcommand{\bfmu}{{\underline{\mu}}}
\newcommand{\bfnu}{{\underline{\nu}}}
\newcommand{\bfV}{{\underline{V}}}

\begin{document}
\mbox{}
\vspace{1cm}

\noindent
\begin{center}
  {\LARGE On\\
   Quantum Cohomology Rings\\
   of\\
   Fano Manifolds\\
   and\\[12pt]
   a Formula of Vafa and Intriligator}
\end{center}
\vspace{1.5cm}

\noindent
\begin{center}
  {\large Bernd Siebert}\footnote{
   On leave from: Mathematisches Institut, Bunsenstr.\ 3--5,
   37073 G\"ottingen, Germany}
   \hspace{1cm}and\hspace{1cm}{\large Gang Tian}\\
   Courant Institute of Mathematical Sciences\\
   251 Mercer Street, New York, NY 10012, USA\\[5pt]
   E-mail: siebert@cims.nyu.edu and tiang@cims.nyu.edu
\end{center}
\vspace{2.5cm}


\noindent
{\Large\bf Introduction}\\[2ex]
Quantum multiplications on the cohomology of symplectic manifolds were
first proposed by the physicist Vafa \cite{vafa} based on Witten's
topological sigma models \cite{witten1}. In \cite {ruan/tian}, Ruan and
the second named author gave a mathematical construction of quantum
multiplications on cohomology groups of {\em positive} symplectic
manifolds  (cf.\ Chapter~1). The construction uses the
Gromov-Ruan-Witten invariants (GRW-invariants in the sequel) for
semi-positive symplectic manifolds previously defined by Ruan
\cite{ruan}. A large class of such manifolds is provided by Fano
manifolds (complex manifolds with ample anti-canonical bundle),
including e.g.\ complex projective spaces, Del Pezzo surfaces and
Grassmann manifolds. Let $M$ be a Fano (or positive symplectic)
manifold. The quantum cohomology $\qrM$ is just the cohomology space
$H^*(M,\cz)$ with a (non-graded) associative, anti-commutative
multiplication, the quantum multiplication.  This multiplication
depends on the choice of a (complexified) K\"ahler class $[\omega]$ on
$M$. Its homogeneous part (the ``weak coupling limit''
$\lambda\cdot[\omega]$, $\lambda\rightarrow\infty$) is the usual cup
product.

In this note we observe that quantum cohomology rings have a nice
description in terms of generators and relations: If
$H^*(M,\cz)=\cz[X_1,\ldots,X_N]/(f_1,\ldots,f_k)$ is a presentation of
the cohomology ring (for simplicity we assume $\deg X_i$ even for the
moment) then $\qrM=\cz[X_1,\ldots,X_N]/(f_1^{[\omega]},
\ldots,f_k^{[\omega]})$, where $f_1^{[\omega]},\ldots,f_k^{[\omega]}$
are just the polynomials $f_1,\ldots,f_k$ evaluated in the quantum
ring associated to $[\omega]$ (Theorem~\ref{presentation}). The
$f_i^{[\omega]}$ actually are real-analytic in $[\omega]$ and thus have
a natural analytic extension to $H^{1,1}(M)$, and the quantum cohomology
rings fit together to form a flat analytic family over $H^{1,1}(M)$.

As an application of this observation we shall compute the
quantum cohomology of the Grassmannians.  The calculations for $\grass$
reduce to the single quantum product $c_k\qwedge s_{n-k}$ of the top
non-vanishing Chern  respectively Segre class of the tautological
$k$-bundle.  This turns out to be a simple exercise in linear algebra.
We derive
\begin{theorem}\label{qrGrass}
  Let $S$ be the tautological bundle over the Grassmann manifold
  $\grass$ of complex $k$-planes in $\cz^n$, $[\omega]=-\lambda\cdot
  c_1(S)$ a $(1,1)$-class on $\grass$ ($\lambda\in\rz_{>0}\
  \Leftrightarrow\ [\omega]$ K\"ahler).  Then
  \[
    \qrG=\cz[X_1,\ldots,X_k]/(Y_{n-k+1},\ldots,Y_{n-1},Y_n+(-1)^k
         e^{-\lambda}),
  \]
  where $X_i$ corresponds to the $i$-th Chern class $c_i(S)$ and
  the $Y_j$ (corresponding to the $j$-th Segre class of $S$) are
  given recursively by ($X_i=0$ for $i>k$)
  \[
    Y_j=-Y_{j-1}\cdot X_1-\ldots-Y_1\cdot X_{j-1}-X_j.
  \]
\end{theorem}
This is in fact the form previously derived by Vafa using arguments
from Quantum Field Theory \cite{vafa0}, \cite{vafa} (but note the sign
of the quantum contribution). It also reduces the genus zero case of a
nice conjectural formula for higher GRW-invariants of $\grass$ stated
in \cite{vafa0} and considered in a broader context in \cite{intr} to a
purely algebraic problem.  A mathematical formulation of the conjecture
in terms of intersection theory on certain algebraic compactifications
of moduli spaces of maps from a (fixed) Riemann surface to $\grass$ was
given by Bertram, Daskalopoulos and Wentworth \cite{bertram.etal}.
They were able to verify the formula for genus one and $k=2$.
\vspace{3ex}

This was the situation known to the authors around New Year 93/94.  A
preliminary version of this paper was spread in a limited number in the
second half of January 94. Since then some progress on related problems
came to our knowledge: Witten (in a paper finished already in mid
December 93) computed the quantum product $c_k\wedge s_{n-k}$,
up to a check of genericity conditions. His arguments are more or less
identical with our Lemma~\ref{holinvariant} (the content of
Lemma~\ref{minratcurve} was taken for granted by him and should in fact
be known classically). The underlying algebraic situation (the content
of Chapter~2 and Chapter~4 of the present paper), however, and thus a
proof of the structure of the quantum cohomology of $\grass$ (even to a
physical level of rigor), was left open.  What nevertheless was
inspiring for us, was his use of residues in the interpretation of
coefficients of quantum products. The analysis of these residues lead
to the first complete proof of the intriguing formula of Vafa and
Intriligator for higher GRW-invariants of $\grass$
(Theorem~\ref{VafIntForm}). On the way we found a nice (new?)
interpretation of (Severi-Grothendieck-Griffiths) residues when there
are ``no components at infinity'' (Proposition~\ref{residues}). An
analogue of the formula of Vafa and Intriligator thus holds for all
Fano manifolds whose cohomology ring is a commutative complete
intersection, including the case where the corresponding sigma model
has a Landau-Ginzburg description (Theorem~\ref{main4}). Later we
realized that certain arguments in Chapter~4 are (encoded in physical
language) already contained in the derivations of \cite{vafa0}.
Moreover, some of the considerations in Chapter~2 have been found
independently by Piunikhin \cite{piuni}.

Further partial support for the formula of Vafa/Intriligator were
obtained in \cite{Rafietal}. The quantum cohomology of (absolute) flag
manifolds (containing $\grass$ as special case) has been calculated in
\cite{GiventalKim} (complete flags) and \cite{Astash} (general case)
{\em assuming} the existence of an equivariant version of quantum
cohomology with certain functorial properties. After completion of this
paper we received a preprint of Bertram containing the reduction to lower
genus for the case of $\grass$ \cite{bertram}.
\vspace{3ex}

The first named author gratefully acknowledges support from the DFG
(enabling him to visit the Courant Institute during the academic year
93/94) and he wishes to thank the Courant Institute for hospitality.
The second author is partly supported by a NSF grant
and an Alfred P.\ Sloan Fellowship.


\section{Definition of quantum multiplications}
In this section, we recall the definition of quantum multiplications
given in \cite{ruan/tian}. The definition uses the
GRW-invariants as defined in \cite{ruan}.

A symplectic manifold $(M,\omega)$ of dimension $2n$ is called {\em
(semi)-positive} if for any $R=f_*[S^2]\in H_2(M,\gz)$,
$f:S^2\rightarrow M$, with $[\omega](R)>0$ either $c_1(M)(R)>0$
($c_1(M)(R)\ge0$). Any Fano manifold is positive w.r.t.\ any K\"ahler
form. Let $J$ be a generic almost complex structure on $M$ tamed by
$\omega$ (i.e.\ $\omega(X,JX)>0$ $\forall\, X\in TM\setminus\{0\}$).  A
$J$-holomorphic curve is a smooth map $f:\Sigma\rightarrow M$
satisfying $J\circ Df = Df\circ j$, where $\Sigma$ is a Riemann surface
(of genus $g$, say) and $j$ is its standard complex structure. This
last equation is a Cauchy-Riemann equation $\dbar f=0$,
$\dbar=\frac{1}{2}(D-J\circ D\circ j)$. For a technical reason it is
convenient to look at the inhomogeneous equation $\dbar f=\gamma$ with
$\gamma$ (the pull-back to $\Sigma$ of) a section of an appropriate
bundle over $\Sigma\times M$.  Solutions of this equation are called
perturbed or $(J,\gamma)$-holomorphic curves. The GRW-invariant can be
defined as follows:

Let $R\in H_2(M,\gz)$ and $B_1,B_2,\ldots,B_s$ be integral
homology classes in $H_*(M,\gz)$ satisfying:
\[
  \hspace{3cm}\sum_{i=1}^s (2n-\deg B_i)= 2c_1(M)(R)+2n(1-g).
  \hspace{3cm}\mbox{(dim)}
\]
Every integral homology class can be represented by a so-called
pseudo-manifold (a certain simplicial complex, cf.\ \cite{ruan/tian}).
For simplicity, we shall also use $B_i$ to denote the pseudo-manifolds
representing these homology classes. Then if $(J,\gamma)$ is generic
and $B_i$ are in sufficiently general position (transversal w.r.t.\ the
evaluation map $\Sigma\times\{f\} \rightarrow M$), for $s$ generic
points $t_1,\ldots,t_s \in\Sigma$ there are only finitely many
$(J,\gamma)$-holomorphic curves $f$ satisfying: $f(t_i)\in B_i$,
$i=1,\ldots,s$ and $f_*[\Sigma]= R$. We define $\tilde\Phi
_{(R,\omega)}(B_1,\ldots,B_s)$ to be the algebraic sum of such $f$ with
appropriate sign according to the orientation. One can prove that
$\tilde\Phi_{(R,\omega)} (B_1,\ldots,B_s)$ is independent of the
choices of $J$, $j$, $\gamma$, and pseudo-manifolds representing
$B_1,\ldots,B_s$ provided the $B_i$ are transversal to the Gromov
boundary of the compactified moduli space of $(J,\gamma)$-holomorphic
curves. The Gromov boundary here consists of curves $C=f(\Sigma)\cup
g_1(S^2)\cup\ldots\cup g_k(S^2)$ with:  $C$ connected; $f$ a
$(J,\gamma)$-holomorphic curve; $g_\nu$ a $(J,0)$-holomorphic rational
curve (a ``bubble''); $R=f_*[\Sigma]+\sum_\nu a_\nu {g_\nu}_*[S^2]$ for
some $a_\nu\in\nz$. Transversality means that there are no such
curves intersecting each $B_i$. Furthermore, $\tilde\Phi_{(R,\omega)}
(B_1,\ldots,B_s)$ depends only on the deformation class $\{\omega\}$ of
$\omega$. Therefore we obtain the GRW-invariant
$\tilde\Phi_{(R,\{\omega\})}(B_1,\ldots,B_s)$.  As a matter of notation
we define $\tilde\Phi_{(R,\omega)}(B_1,\ldots,B_s)$ to be zero unless
the dimensions match (dim).

We now make a simple but important remark:  Let $J$ be an almost
complex structure such that any $J$-holomorphic curve with
$f_*[\Sigma] = R$ is regular, i.e.\ the linearization of the
Cauchy-Riemann operator at $f$ has trivial cokernel. Then $(J,0)$ is
generic and we can use exact $J$-holomorphic curves to compute the
invariant $\tilde\Phi_{(R, \{\omega\})}(B_1,\ldots,B_s)$. In case $J$ is
integrable and $f$ is an immersion, regularity at $f$ is equivalent to
the vanishing of $H^1(C,N_{C|M})$, $C=f(\Sigma)$ \cite[2.1.B.]{gromov}. In
particular the invariant is easily computable in certain cases as
follows:
\begin{lemma}\label{genericcrit}
  Let $(M,\omega)$ be a K\"ahler manifold and $R\in H_2(M,\gz)$ s.th.\
  any holomorphic curve $C\subset M$ homologous to $R$ is non-singular
  and has $H^1(C,N_{C|X})=0$. Let $B_1,\ldots,B_s\subset M$ be complex
  submanifolds transversal to the evaluation map and the Gromov boundary.
  Then
  \[
    \tilde\Phi_{(R,\{\omega\})}(B_1,\ldots,B_s)=\sum_C\,
    \sharp(B_1\cap C)\cdot\sharp(B_2\cap C)\cdot\ldots\cdot\sharp(B_s\cap C),
  \]
  where the sum is over all holomorphic curves $C$ homologous to $R$.
\end{lemma}
Note that $\sharp(B_1\cap C)\cdot\ldots\cdot\sharp(B_s\cap C)$ is
precisely the number of ways to parametrize $C$ by $f:\Sigma\rightarrow
C$ with $t_i$ mapping to $B_i$, $i=1,\ldots,s$. Moreover, there are no
signs occurring in the formula since all spaces involved are
canonically oriented by their complex structure.

\begin{rem}\label{algmod}\rm
(This will only be used to comment on the relationship of our work with
\cite{bertram.etal}.)
More generally, if the moduli space $\calm_{R,\Sigma}$ (for shortness
$\calm$ in the sequel) of holomorphic maps
$f:\Sigma\rightarrow M$ with $f_*[\Sigma]=R$ has the expected dimension
$c_1(M)(R)+n(1-g)$ (equivalently, $H^1(\Sigma,f^*T_M)=0$ for almost
all $f$) the GRW-invariants $\tilde\Phi_R$ have an interpretation as
intersection products of certain compactifications of $\calm$: $M$ being
projective algebraic
there clearly exist projective compactifications of $\calm$. Choose one,
say $\overline\calm$. Then the evaluation map
\[
  {\rm ev}:\Sigma\times\calm\longrightarrow M,\ \ \
  (t,f)\longmapsto f(t)
\]
extends rationally to $\overline\calm$. Hence there is a modification
$\pi:\Gamma\rightarrow\Sigma\times\overline\calm$, biregular over
$\Sigma\times\calm$ such that ${\rm ev}$ has an extension $\tilde{\rm ev}:
\Gamma\rightarrow M$. For generic $t\in\Sigma$, $\Gamma_t:=\pi^{-1}(\{t\}\times
\calm)$ is a modification of $\overline\calm$. We choose a common
desingularization
$\widehat\calm$ of all the $\Gamma_t$, $t\in\Sigma$ generic (blow-up the sum
of ideals $\cali_t$ whose blow-up desingularizes $\Gamma_t$, then
desingularize).
Let $\iota_t:\widehat\calm\rightarrow\Gamma_t \hookrightarrow\Gamma$. Then
for $\beta_i\in H^*(M,\gz)$ Poincar\'e-dual to $B_i$ and generic choices of
$t_1,\ldots,t_s\in\Sigma$
\[
  \tilde\Phi_R(B_1,\ldots,B_s)
  =(\iota_{t_1}^*\Phi^*\beta_1\wedge\ldots\wedge\iota_{t_s}^*\Phi^*\beta_s)
  [\widehat\calm].
\]
In fact, choosing $B_1,\ldots,B_s$ as generic pseudo-manifolds then
$(\Phi\circ \iota_{t_i})^{-1}(B_i)$ is Poincar\'e-dual to
$\iota_{t_i}^*\Phi^*\beta_i$ and the bordism argument of \cite{ruan}
should generalize to show
\[
  \tilde\Phi_R(B_1,\ldots,B_s)=\sharp(\Phi\circ\iota_{t_1})^{-1}(B_1)\cap
  \ldots\cap(\Phi\circ\iota_{t_s})^{-1}(B_s),
\]
($\sharp$ means algebraic sum with appropriate signs keeping track of
the orientations). A proper proof has to provide a bordism to a generic
situation carefully dealing with the singularities of $\calm$ as well
as with transversality and orientations.  Details will be given in
\cite{ruan/tian}.
\qed
\end{rem}

To define a quantum multiplication on $H^*(M,\gz)$, we introduce
the real-valued invariant
\[
  \tilde{\Phi}_{[\omega]}(B_1,\ldots,B_s)=\sum_{R\in H_2(M,\gz)}
  \tilde{\Phi}_{(R,\{\omega\})}(B_1,\ldots,B_s)\,e^{-[\omega](R)}.
\]
In general there might be infinitely many terms contributing to the sum
(e.g.\ in the important and interesting case of Calabi-Yau manifolds)
and one faces a non-trivial convergence problem.  In the positive case,
however, this sum is actually finite due to the dimension condition
(dim). Let us assume $M$ positive for the following (but see
Remark~\ref{cy-rem}).  Note also that letting $\omega$ vary one gets
back all the invariants $\tilde{\Phi}_{(R,\{\omega\})}(B_1,\ldots,B_s)$
by solving some system of linear equations.  The quantum multiplication
$\qwedge$ on $H^*(M,\rz)$ is then characterized by the equation
\[
  (\alpha\qwedge\beta)[A]=\tilde{\Phi}_{[\omega]}(\alpha^\vee,
  \beta^\vee,A),
\]
where $\alpha$, $\beta\in H^*(M,\gz)$ and $\ ^\vee$ means Poincar\'e-dual.
In terms of a basis $\{A_i\}$ for the torsion free part of
$H_*(M,\gz)$ and $\{\alpha_i\}$ the Poincar\'e-dual basis of
$H^*(M,\gz)$ we may state this more explicitely as
\[
  \alpha_i\qwedge\alpha_j=\sum_{k,l}\eta^{lk}\skp\tilde{\Phi}_{[\omega]}(A_i,
  A_j,A_k)\,\alpha_l,
\]
with $(\eta^{lk})_{kl}$ inverse to the intersection matrix
$(\eta_{ij})_{ij}=(A_i\cdot A_j)_{ij}$ ($\sum_k
\eta_{ij}\eta^{kj}=\delta_{ik}$). The associativity of the quantum
multiplication is highly non-trivial and shown by a careful analysis
of the degeneration of rational curves (following an idea of Witten)
in \cite{ruan/tian}.

A trivial, but decisive feature of this definition, that we are going
to exploit, is that its homogeneous part reduces to the cup product
\[
  \alpha_i\wedge\alpha_j=\sum_{k,l}\eta^{lk}\cdot(A_i\cdot A_j\cdot A_k)
  \,\alpha_l.
\]


\section{A presentation for quantum cohomology}
Denote by $\cz\langle X_1,\ldots,X_n\rangle$ the graded anticommutative
$\cz$-algebra with generators $X_i$ of degree $d_i$, i.e.\ with $X_i
X_j=(-1)^{d_i\cdot d_j}X_j X_i$. If $m$ of the $X_i$ have odd degree
this is isomorphic to $\left(\Lambda^*\cz^m\right) \otimes\left(
\mbox{Sym}^*\cz^{n-m}\right)$.  We will call elements of this algebra
ordered polynomials (since in addition to the coefficients one has to
select an order among the factors in a monomial to determine its
sign).

Let $(M,\omega)$ be a positive symplectic manifold and
\[
  H^*(M,\cz)=\cz\langle X_1,\ldots,X_n\rangle/(f_1,\ldots,f_k)
\]
be a presentation of the cohomology ring, $f_i=\sum_{|J|=\deg f_i}
a_{iJ}X^J$. (We use multiindex notation $J=(j_1,\ldots,j_n)$,
$X^J=X_1^{j_1}\wedge\ldots\wedge X_n^{j_n}$, $|J|=\sum_{i=1}^n j_i d_i$
etc.) Denote by $\qwedge$ the product in the quantum cohomology $\qrM$.
To distinguish clearly between calculations in $H^*(M,\cz)$ and in
$\qrM$, we use a $\widehat{\ }$ to mark elements of the quantum
cohomology ($\widehat{\ }$ might be thought of as a $\cz$-linear map
$\cz\langle X_1,\ldots,X_n\rangle\rightarrow\qrM$).  An ordered polynomial
$f$ may be evaluated on generators of the quantum cohomology provided
the degrees match.  We write $f(\hat X_1,\ldots,\hat X_n)$ as in the
commutative case.  Especially, $\hat X^J:=\hat
X_1\qwedge\ldots\qwedge\hat X_1\qwedge\ldots\qwedge\hat X_n
\qwedge\ldots\qwedge\hat X_n$ with $\hat X_\nu$ occurring $j_\nu$-often
(note the difference from \mbox{$(X^J)\hat{\,}$ ).}
\begin{lemma}\label{generate}
  $\hat X_1,\ldots,\hat X_n$ generate $\qrM$.
\end{lemma}
\pf
By induction on the degree. So assume $\hat X_1,\ldots,\hat X_n$
generate $\qrM$ up to degree $d-1$.  We want to show that the
monomials $(X^J)\hat{\,}$, $|J|=d$, can be written as linear combination
of (quantum) products of the $\hat X_\nu$. But by definition of the
multiplication in $\qrM$
\[
  \hat X^J= (X^J)\hat{\,}+\sum_{|I|<d}a_I (X^I)\hat{\,}\, ,
\]
and by induction hypothesis $(X^I)\hat{\,}=\sum_{K\le I}b_{IK}\hat X^K$, so
\[
  (X^J)\hat{\,}=\hat X^J-\sum_{K\le I\atop |I|<d}a_I b_{IK}\hat X^K.
\]
\vspace{-1cm}

\qed
\vspace{2ex}

Next we are calculating $f_i(\hat X_1,\ldots,\hat X_n)$ (i.e. in the
quantum cohomology).  Since $f_i$ is a relation in the cohomology ring, the
($\deg f_i$)-term vanishes. By the lemma we find an ordered polynomial
$g_i^{[\omega]}$ (depending on $[\omega]$) in $n$ variables
(say $T_1,\ldots,T_n$, $\deg T_i=\deg X_i$) of {\em lower degree} with
\[
  f_i(\hat X_1,\ldots,\hat X_n)=g_i^{[\omega]}(\hat X_1,\ldots,\hat X_n)
  \hspace{1cm}\mbox{in }\qrM.
\]
Thus
$f_i^{[\omega]}(T_1,\ldots,T_n):=f_i(T_1,\ldots,T_n)-
g_i^{[\omega]}(T_1,\ldots,T_n) \in\cz\langle T_1,\ldots,T_n\rangle$ is
a non-trivial relation between $\hat X_1,\ldots,\hat X_n$.
\begin{theorem}\label{presentation}
  $\qrM=\cz\langle T_1,\ldots,T_n\rangle/(f_1^{[\omega]},
  \ldots,f_k^{[\omega]})$.
\end{theorem}
\pf
Let $\calj\subset\cz\langle T_1,\ldots,T_n\rangle$ be the ideal of
relations between the $\hat X_1,\ldots,\hat X_n$. Then
$(f_1^{[\omega]},\ldots,f_k^{[\omega]})\subset \calj$ and $\qrM=\cz\langle
T_1,\ldots,T_n\rangle/\calj$ by the lemma. Let $F\in \calj\setminus\{0\}$.
Expand $F=F_d+F'$ with $F_d\neq0$ (weighted) homogeneous of degree $d$,
$d>0$, and $\deg F'<d$. Then
\[
  F_d(\hat X_1,\ldots,\hat X_n)=-F'(\hat X_1,\ldots,\hat X_n)
\]
in the quantum cohomology for $F$ is a relation. But the highest degree
($\deg=d$) contribution to $F_d(\hat X_1,\ldots,\hat X_n)$ is just
$(F_d(X_1,\ldots,X_n))\hat{\,}\,$. Its vanishing implies
$F_d\in(f_1,\ldots,f_k)$, i.e.\ $F_d=\ph(f_1,\ldots,f_k)$, $\ph$ a
polynomial in $k$ variables. Thus
\[
  \ph(f_1^{[\omega]},\ldots,f_k^{[\omega]})=F_d+F''
  \hspace{1cm}\mbox{in }\cz\langle T_1,\ldots,T_n\rangle
\]
with $\deg F''<d$, and we may write
$F=\ph(f_1^{[\omega]},\ldots,f_k^{[\omega]}) +F'-F''$, $\deg F'-F''<
d=\deg F$. Proceeding by induction on the degree we finally see
$F\in(f_1^{[\omega]},\ldots,f_k^{[\omega]})$,
i.e.\ $\calj=(f_1^{[\omega]},\ldots,f_k^{[\omega]})$ as claimed.
\qed
\vspace{2ex}
\begin{rem}\label{cy-rem}\rm
In the positive case $c_1(M)>0$ the contributions with $R\neq0$ give rise
to terms of lower degree by the dimension condition, so we were able to fix
$[\omega]$ and argue by considerations on the degree. Modulo
convergence problems in the semi-positive case mentioned at the end
of the previous chapter, Theorem~\ref{presentation} however remains valid.
Arguments involving the degree can be replaced by the linear
independence of terms $e^{-\lambda\cdot t}$ for various $\lambda$
in the algebra $\cz[[t]]$.

Another approach, which works without any positivity condition, is to
use formal power series over $H_2(M,\gz)$ (cf.\ also \cite{piuni}): Let $N$
be the (obvious) completion of the group ring over $H_2(M,\gz)$. Then
$H^*(M,\gz)\otimes_\gz N$ appears as natural domain of definition of
quantum products for more general manifolds. Our arguments provide a
presentation of the quantum cohomology ring as quotient of $\gz\langle
T_1,\ldots,T_n\rangle\otimes_\gz N$ by ordered polynomials with coefficients
in $N$. Note also that in the positive case $\qrM$ may be viewed as
homeomorphic image of this formal quantum cohomology ring by sending
$R\in H_2(M,\gz)$ to $e^{-[\omega](R)}$.
\qed
\end{rem}


\section{Grassmann manifolds}
As is well-known (e.g.\ \cite[Ex.~14.6.6]{fulton}) the cohomology ring
(in fact even the Chow ring) of the Grassmann variety $\grass$ of
$k$-planes in $\cz^n$ has a presentation
\[
  H^*(\grass,\cz)=\cz[c_1,\ldots,c_k]/(s_{n-k+1},\ldots,s_n),
\]
with $c_i$ corresponding to the Chern classes of the tautological
$k$-bundle $S$ and $s_j$ to its Segre classes viewed as polynomials in
$c_1,\ldots,c_k$ via
\begin{eqnarray*}
  &(1+c_1+\ldots+c_k)(1+s_1+s_2+\ldots)=1,&\\
  &\mbox{i.e. } s_j=-s_{j-1}c_1-\ldots-s_1c_{j-1}-c_j.&\hspace{3cm}(*)
\end{eqnarray*}
Note that $s_j$ is also the $j$-th Chern class of the universal
quotient bundle $Q$, which has rank $n-k$, so $s_j=0$ in
$H^*(\grass,\cz)$ for $j>n-k$. In fact, as a polynomial in $c_i$, $s_j$
lies in the relation ideal for $j>n$ by the recursion formula ($*$).
Under the canonical isomorphism $\Phi:\grass\simeq{\rm G}(n-k,n),\
\Lambda\mapsto(\cz^n/\Lambda)^*$, $S$ corresponds to the dual of the
universal quotient bundle $Q'$ on ${\rm G}(n-k,n)$ and $Q$ to the dual of the
tautological $(n-k)$-bundle $S'$.  Thus $c_i$ and $s_j$ exchange (up to sign)
their roles and one might as well write
\[
  H^*(\grass,\cz)=\cz[s_1,\ldots,s_{n-k}]/(c_{k+1},\ldots,c_n),
\]
this time with $c_i$ polynomials in $s_1,\ldots,s_{n-k}$ (this
presentation is actually better adapted to Schubert calculus,
i.e.\ geometry, see below). These remarks are made to emphasize the symmetry
between $c_i$ and $s_j$. Note also that the generators all have
even degree, so we need not worry about questions of sign.

Schubert calculus (e.g.\ \cite[\S~14.7]{fulton}) provides a basis of
$H^*(\grass,\cz)$ as a $\cz$-vector space (indeed a basis of integral
homology/cohomology as $\gz$-module), indexed by tuples
$(\lambda_1,\ldots,\lambda_k)$,
$n-k\ge\lambda_1\ge\ldots\ge\lambda_k\ge0$, via
\[
  \{\lambda_1,\ldots,\lambda_k\}:=\det(s_{\lambda_i+j-i})_{1\le i,j\le k},
\]
i.e.\ by evaluating the Schur polynomial (S-function) associated to
$(\lambda_1,\ldots,\lambda_k)$ on the Segre classes of the tautological
$k$-bundle ($s_j=0$ for $j\not\in\{0,\ldots,n-k\}$).
$\{\lambda_1,\ldots,\lambda_k\}$ is (weighted) homogeneous of degree
$2\sum\lambda_i$.  The Chern  respectively Segre classes are given by
\begin{eqnarray*}
  c_i&=&(-1)^i\{1,\ldots,1,0,\ldots,0\}\hspace{3ex}\mbox{(``1'' $i$-times)},\\
  s_j&=&\{j,0,\ldots,0\}
\end{eqnarray*}
(for $c_i$ see \cite[Lemma~14.5.1]{fulton} or do an easy induction).
The connection to classical Schubert calculus is given by
Poincar\'e-duality. In fact, $\{\lambda_1,\ldots,\lambda_k\}\cap[\grass]$
may be represented by the Schubert varieties
\[
  \Omega_\bfV(n-k+1-\lambda_1,\ldots,n-k+i-\lambda_i,\ldots,n-\lambda_k).
\]
To define the latter one has to fix a flag $\bfV=(V_1,\ldots,V_n)$,
$V_1\subset\ldots\subset V_n=\cz^n$ of linear subspaces of
$\cz^n$, $\dim V_i=i$. Then for $(a_1,\ldots,a_k)$,
$0\le a_1<\ldots<a_k\le n$
\[
  \Omega_\bfV(a_1,\ldots,a_k):=\{\Lambda\in \grass\ |
  \dim\Lambda\cap V_{a_i}\ge i, 1\le i\le k\}.
\]
The homology {\em class}
$(a_1,\ldots,a_n):=[\Omega_\bfV(a_1,\ldots,a_k)]$ is independent of
the choice of flag. We write $\{\bflm\}^\vee:=(n-k+1-\lambda_1,\ldots,
n-k+i-\lambda_i,\ldots,n-\lambda_k)\in H_*(\grass,\cz)$. What is
important for us is that $[\grass]$ is Poincar\'e-dual to
$\{0,\ldots,0\}=1$ (trivial) and that the class of a point $[*]$ is
Poincar\'e-dual to $\{n-k,\ldots,n-k\}=s_{n-k}^k=(-1)^{n-k}c_k^{n-k}$
(which makes sense as $\dim \grass=k(n-k)$, but is less trivial: What
is the Poincar\'e-dual of $c_1^{k(n-k)}$?). The intersection of classes of
complementary dimension is especially easy:  $|\bflm|+|\bfmu|=k(n-k)$,
then (``duality theorem''):
\[
  \left(\{\bflm\}\wedge\{\bfmu\}\right)[G]
  =\{\bflm\}\cap\{\bfmu\}^\vee
  =\{\bflm\}^\vee\cdot\{\bfmu\}^\vee
  =\left\{\begin{array}{lcl}1&,&\mbox{ if }\bfmu=\bflm^*.\\
                           0&,&\mbox{ otherwise}.\end{array}\right.
\]
with $\bflm^*:=(n-k-\lambda_k,\ldots,n-k-\lambda_1)$. In particular, we
get $\eta_{\bflm\bfmu}=\delta_{\bflm\bfmu}$ for the intersection matrix.
This ends our collection of facts concerning Grassmannians.
\vspace{2ex}

The product in $\qrG$ now reads
\[
  \widehat{\{\bflm\}}\qwedge\widehat{\{\bfmu\}}=\sum_{R\in H_2(\grass,
  \gz)}\sum_{\bfnu}\tilde\Phi_R\left(\{\bflm\}^\vee,
  \{\bfmu\}^\vee,\{\bfnu\}^\vee\right)\widehat{\{\bfnu^*\}}
  e^{-[\omega](R)}.\hspace{1cm}{(**)}
\]
Next we observe that $H^{1,1}(\grass)$ is spanned by the single
class $\{1,0,\ldots,0\}=s_1=-c_1$ and that $\det(Q)=(\det
S)^*$ is actually very ample, the corresponding embedding being the
Pl\"ucker embedding $\iota:\grass\hookrightarrow \pr(\Lambda^k\cz^n)$
(so $s_1=\iota^*c_1\left(\calo_{\pr(\Lambda^k\scz^n)}(1)\right)$).
Dually, $H_2(\grass,\gz)$ is spanned by the single class
$\{n-k,\ldots,n-k,n-k-1\}^\vee=(1,\ldots,k-1,k+1)=:[L]$, and the sum
over $H_2(\grass,\gz)$ is just running through the set
$\{d\cdot[L]\mid d\in {\nz}_0\}$. Since
$T_{\grass}\simeq\Hom(S,Q)\simeq S^*\otimes Q$, as one easily verifies
by using the standard local coordinates on $\grass$,
$c_1(\grass)=\rk(Q)\cdot c_1(S^*)+\rk(S^*)\cdot c_1(Q)=n\cdot s_1$.

Recall from Chapter~1 that $\tilde\Phi_{d\cdot[L]}\left(\{\bflm\}^\vee,
\{\bfmu\}^\vee,\{\bfnu\}^\vee\right)=0$ unless
$|\bflm|+|\bfmu|+|\bfnu|=k(n-k)+c_1(\grass)
\left(d\cdot[L]\right)=k(n-k)+d\cdot n$. In calculating
$s_{n-k+i}(\hat{c}_1,\ldots,\hat{c}_k)$ (i.e.\ in $\qrG$,
$\hat{c}_i$ the generators corresponding to $c_i$ as in Chapter~2) we
have $|\bflm|+|\bfmu|\le n$, and ``='' only in the case of $s_n$.  But
$|\bfnu|\le k(n-k)=\dim \grass$, so we get no quantum contributions
besides in case of $s_n$ with $|\bfnu|=k(n-k)$ and $d=1$. Thus
letting $\underline{\hat c}=(\hat{c}_1,\ldots,\hat{c}_k)$,
$s_{n-k+1}(\underline{\hat c})=\ldots=s_{n-1}(\underline{\hat c})=0$ in
$\qrG$ and
\begin{eqnarray*}
  s_n(\underline{\hat c})
  &=&-\hat{c}_1 s_{n-1}(\underline{\hat c})
  -\ldots-\hat{c}_{k-1}s_{n-k+1}(\underline{\hat c})
  -\hat{c}_k s_{n-k}(\underline{\hat c})\\
  &=&-\hat{c}_k s_{n-k}(\underline{\hat c})\ =\
  -\tilde\Phi_{[L]}\left(c_k^\vee,s_{n-k}^\vee,[*]\right)
  \cdot e^{-[\omega](L)},
\end{eqnarray*}
where for the last equality we have used $(**)$.
Taking into account Theorem~\ref{presentation} to prove
Theorem~\ref{qrGrass} we are left with
\begin{prop}\label{Phi_L}
  $\tilde\Phi_{[L]}\left(c_k^\vee,s_{n-k}^\vee,[*]\right)=(-1)^k$.
\end{prop}
For the proposition we are first classifying holomorphic curves homologous
to $[L]=(1,\ldots, k-1,k+1)$.
\begin{lemma}\label{minratcurve}
  Let $C\subset \grass$ be a (rational) curve homologous to
  $(1,\ldots,k-1,k+1)$. Then $C$ is a Schubert variety
  $\Omega_{\bfV}(1,\ldots,k-1,k+1)$, i.e.\ there are linear subspaces
  $U\subset W\subset\cz^n$, $\dim U=k-1$, $\dim W=k+1$, with
  $C=\{\Lambda\in \grass\mid U\subset\Lambda\subset W\}$.
\end{lemma}
\pf
Since $s_1[C]=1$ by the duality theorem, $\deg\iota(C)=1$, so the image
of $C$ under the Pl\"ucker embedding
$\iota:\grass\rightarrow\pr(\Lambda^k\cz^n)$ is a linear $\pr^1$. The
image of $\grass$ consists precisely of (rays of) {\em decomposable}
vectors $v_1\wedge\ldots\wedge v_k=\iota(\langle v_1,\ldots,v_k\rangle)
\in\Lambda^k\cz^n$. Locally the vectors $v_1,\ldots,v_k$ may be chosen
to vary smoothly with $\Lambda\in \grass$: In fact, the standard
(affine) coordinate neighbourhood of $\Lambda\in \grass$ is
$\Hom(\Lambda,\cz^n/\Lambda)$ with
$\Psi:\Hom(\Lambda,\cz^n/\Lambda)\rightarrow \grass$,
$\ph\mapsto\langle v+\ph(v) \mid v\in\Lambda\rangle$ (in particular
$0\in\Hom(\Lambda,\cz^n/\Lambda)$ corresponds to $\Lambda$). Now fixing
a basis $v_1,\ldots,v_k$ of $\Lambda$, $\iota\circ\Psi$ may be
represented (lifted to $\Lambda^k\cz^n$) by
\[
  \ph\mapsto\left(v_1+\ph(v_1))\wedge\ldots\wedge(v_k+\ph(v_k)\right).
\]
Thus choosing $\Lambda'\in C$ sufficiently close to $\Lambda\in C$ and letting
$e_1,\ldots,e_l$ be a basis of $\Lambda\cap\Lambda'$, completed by $e_{l+1},
\ldots,e_k$ and $e'_{l+1},\ldots,e'_k$ to a basis of $\Lambda$ and
$\Lambda'$ respectively, we have for $t\in\cz$ small
\[
  e_1\wedge\ldots\wedge e_l\wedge e_{l+1}\wedge\ldots\wedge e_k
  +t\cdot e_1\wedge\ldots\wedge e_l\wedge e'_{l+1}\wedge\ldots\wedge e'_k
  =v_1(t)\wedge\ldots\wedge v_k(t)
\]
with $v_i(0)=e_i$ by construction. Taking
$\displaystyle\left.\frac{d}{dt}\right|_{
t=0}$ yields
\[
  e_1\wedge\ldots\wedge e_l\wedge e'_{l+1}\wedge\ldots\wedge e'_k
  =\dot{v}_1(0)\wedge e_2\wedge\ldots\wedge e_k+\ldots+
  e_1\wedge\ldots\wedge e_{k-1}\wedge\dot{v}_k(0).
\]
As one sees by expanding $\dot{v}_i(0)$ in terms of a basis of $\cz^n$
containing $\{e_1,\ldots,e_k,$ $e'_{l+1},\ldots,e'_k\}$ we may gather the
linearly independent terms with $e_k$ and without $e_k$ to form two
equations. The left-hand side of the equation above belongs
to the latter ($\Lambda\neq\Lambda'$ $\Rightarrow$ $l<k$), so we get
\[
  e_1\wedge\ldots\wedge e_l\wedge e'_{l+1}\wedge\ldots\wedge e'_k
  =e_1\wedge\ldots\wedge e_{k-1}\wedge\left(\dot{v}_k(0)-\lambda\cdot
  e_k\right)
\]
for some $\lambda\in\cz$ (s.th.\ $\dot{v}_k(0)-\lambda\cdot e_k$ lies
in the span of the basis vectors different from $e_k$). By linear
independence of wedge products of a basis of $\cz^n$ this shows
$l=k-1$. In view of the linearity of $\iota(C))$ we conclude
\[
  \iota(C)=\left\{ [t\cdot e_1\wedge\ldots\wedge
  e_{k-1}\wedge e'_k+ u\cdot e_1\wedge\ldots\wedge e_{k-1}\wedge
  e_k]\ \Big|\ [t:u]\in\pr^1\right\},
\]
so $C$ is the Schubert variety
$\Omega_{\bfV}(1,\ldots,k-1,k+1)$ belonging to a flag $\bfV$ with
$V_{k-1}=\langle e_1,\ldots,e_{k-1}\rangle=:U$ and $V_{k+1}=\langle
e_1,\ldots,e_{k-1},e_k,e'_k\rangle=:W$.
\qed

\begin{lemma}\label{holinvariant}
  Let $A_1=\Omega_{\underline{V}^1}(n-k,n-k+1,\ldots,n-1)=\{1,\ldots,1\}^\vee
  =(-1)^k c_k^\vee$, $A_2=\Omega_{\underline{V}^2}(1,n-k+2,\ldots,n)=
  \{n-k,0,\ldots,0\}^\vee=s_{n-k}^\vee$,
  $A_3=\{*\}=\Omega_{\underline{V}^3}(1,\ldots,k)
  =\{n-k,\ldots,n-k\}^\vee$, where $\underline{V}^1$,
  $\underline{V}^2$, $\underline{V}^3$ are three transversal flags
  (i.e.\ $\dim V^1_i\cap V^2_j\cap V^3_k=\max\{0,i+j+k-2n\}$).
  Then there is one and only one rational curve $C$ homologous to $[L]$
  and with $C\cap A_i\neq\emptyset$, $i=1,2,3$. Moreover,
  $C\cdot A_1=C\cdot A_2=C\cdot A_3=1$.
\end{lemma}

The lemma can be proved by doing intersection theory on the flag
manifold ${\rm F}(k-1,k+1;n)$, which parametrizes rational curves of minimal
degree by the preceding lemma, and using the two obvious maps
$\pi:{\rm F}(k-1,k+1;n)\rightarrow \grass$ and $p:{\rm F}(k-1,k+1;n)\rightarrow
\grass$ (calculate $(p_*\pi^*[A_1])\cdot(p_*\pi^*[A_2])
\cdot(p_*\pi^*[A_3])$). This method might be appropriate for more
general $A_1,A_2,A_3$, yet in our case an explicit linear algebra
argument is simpler and even more enlightening.
\vspace{2ex}

\noindent\pf
We nee to find three $k$-planes $\Lambda^1$, $\Lambda^2$, $\Lambda^3$
and subspaces $U,W\subset\cz^n$, $\dim U=k-1$, $\dim W=k+1$ with
\begin{enumerate}
\item
  $U\subset\Lambda^1\cap\Lambda^2\cap\Lambda^3\subset
  \Lambda^1+\Lambda^2+\Lambda^3\subset W$.
\item
  $\dim\Lambda^1\cap V^1_{n-k+(i-1)}\ge i$, $i=1,\ldots,k$.
\item
  $V_1^2\subset\Lambda^2$.
\item
  $\Lambda^3=V^3_k$,
\end{enumerate}
(1) says that $\Lambda^1,\Lambda^2,\Lambda^3$ lie on the rational curve
$C$ defined by $U$ and $W$ according to Lemma~\ref{minratcurve},
whereas (2)--(4) rephrase the conditions $\Lambda_i=A_i\cap C$,
$i=1,2,3$ respectively. We are now using transversality of the flags
$\underline{V}^1$, $\underline{V}^2$, $\underline{V}^3$. From (3), (4)
and (1) we readily deduce $W=V_k^3+V_1^2$, and (2) with $i=k$ shows
$\Lambda^1\subset V^1_{n-1}$, so by (4) and (1) we get $U=V_k^3\cap
V^1_{n-1}$. This choice of $U,W$ implies $\Lambda^1 =W\cap V^1_{n-1}$,
$\Lambda^2=U+V^2_1$, $\Lambda^3=V^3_k$. Conversely these
$\Lambda^1,\Lambda^2,\Lambda^3$ fulfill (1)--(4).  \qed
\vspace{2ex}

\noindent
{\em Proof of Proposition~\ref{Phi_L}.}
In view of the preceding lemma and the criterion for multiplicity one
(Lemma~\ref{genericcrit} with $R=[L]$) the only thing remaining to be
checked is transversality of $A_i$ w.r.t.\ the evaluation map
\[
  \Sigma\times\{f:\Sigma\rightarrow\grass\}\rightarrow\grass.
\]
The moduli space $\{f\}$ of holomorphic curves in question is in our
case isomorphic to the flag manifold ${\rm F}(k-1,k+1;n)
=\{(U,W)\}$ by Lemma~\ref{minratcurve} and \ref{holinvariant}.  But
$\Lambda_i=A_i\cap f(\Sigma)$ clearly varies when $f$ (and hence $(U,W)$)
varies as seen explicitely in the proof of Lemma~\ref{holinvariant}.
Note also that the Gromov boundary is empty in this case. This follows either
from the explicit description of the moduli space as flag manifold
or from the fact that $L\in H_2(\grass,\gz)$ is primitive and thus can not
be represented by any reducible holomorphic curve.
\qed

This finishes the proof of Theorem~\ref{qrGrass}.


\section{Higher invariants}
The main results of \cite{ruan/tian} show how to compute higher
GRW-invariants (i.e.\ with more than three entries or for higher genus
Riemann surfaces) from the genus $0$ three-point functions inductively.
Namely, for $g>0$

\[
  {\tilde\Phi}^g_{[\omega]}(B_1,\ldots,B_s)=\sum_{i,j}\eta^{ij}
  {\tilde\Phi}^{g-1}_{[\omega]}(B_1,\ldots,B_s,A_i,A_j),
\]
$(\eta_{ij})$ the intersection matrix with respect to a basis $\{A_i\}$
of $H^*(M,\cz)$. More invariantly, the right-hand side of course is the
trace with respect to $\eta$ of the bilinear form
\[
  H^*(M,\rz)\times H^*(M,\rz)\longrightarrow\rz,\ \ \ (B',B'')\longmapsto
  {\tilde\Phi}^{g-1}_{[\omega]}(B_1,\ldots,B_s,B',B'').
\]
Secondly, for $g=0$ and $1<r<s-1$ (otherwise trivial) we have the
{\em composition law}
\[
  {\tilde\Phi}^0_{[\omega]}(B_1,\ldots,B_s)=\sum_{i,j}\eta^{ij}
  {\tilde\Phi}^0_{[\omega]}(B_1,\ldots,B_r,A_i)
  {\tilde\Phi}^0_{[\omega]}(A_j,B_{r+1},\ldots,B_s),
\]
a trace with respect to $\eta$ as well (for $s=4$ and $r=2$ this
equation states the associativity of quantum products.) Our goal in
this section is to give a closed formula for higher invariants in terms
of the relations $f_1^{[\omega]},\ldots,f_k^{\omega]}$ of the quantum
cohomology ring. Putting all ${\tilde\Phi}^g_{[\omega]}$ together for
different $s$ we get a $\cz$-linear map
\[
  \langle\ \ \rangle_g:\cz\langle X_1,\ldots X_n\rangle\longrightarrow\cz,
  \ \ X_1^{\nu_1}\ldots X_n^{\nu_n}\longmapsto{\tilde\Phi}^g_{[\omega]}
  (X_1^\vee,\ldots,X_1^\vee,\ldots,X_n^\vee,\ldots,X_n^\vee),
\]
the Poincar\'e-dual $X_i^\vee$ of $X_i$ occurring $\nu_i$-times. But
from the composition law denoting by $\alpha_i$ the Poincar\'e-dual of
$A_i$
\[
  \hat X_1^{\nu_1}\qwedge\ldots\qwedge\hat X_n^{\nu_n}=
  \sum_{ij}\eta^{ij}\,{\tilde\Phi}^0_{[\omega]}(X_1^\vee,\ldots,
  X_1^\vee,\ldots,X_n^\vee,A_i)\,\hat\alpha_j,
\]
so $\langle X_1^{\nu_1}\ldots X_n^{\nu_n}\rangle_0={\tilde\Phi}^0_{[\omega]}
(X_1^\vee,\ldots,X_1^\vee,\ldots,X_n^\vee,[M])$ is nothing but the coefficient
of the class $[\Omega]$ of the normalized volume form in $\hat X_1^{\nu_1}
\qwedge\ldots\qwedge\hat X_n^{\nu_n}$. That is, $\langle\ \ \rangle_0$
decomposes
\[
  \cz\langle X_1,\ldots,X_n\rangle\longrightarrow
  \cz\langle X_1,\ldots,X_n\rangle/(f_1^{[\omega]},\ldots,f_k^{[\omega]})
  \simeq H^*_{[\omega]}(M)\stackrel{{\rm top}}{\longrightarrow}
  H^{2n}(M,\cz),
\]
and $H^{2n}(M,\cz)$ is identified with $\cz$ by sending $[\Omega]$ to
$1$. In case $n=k$ and the generators have even degree, i.e.
$H^*(M,\cz)$ a (commutative) complete intersection ring, one can use
higher dimensional residues to express this map more explicitely (we
refer the reader to \cite{griffharr} and \cite{tsikh} for
the general facts on residues to be used):

\sloppy
Recall that the residue of $F\in\cz[X_1,\ldots,X_k]$ with respect to
a polynomial mapping $g=(g_1,\ldots,g_k):\cz^k\rightarrow\cz^k$ with
$g^{-1}(0)$ finite (or equivalently, $\cz[X_1,\ldots,X_k]/(g_1,\ldots,g_k)$
is artinian, i.e.\ finite dimensional as vector space over $\cz$) in
$a\in g^{-1}(0)$ is defined by
\[
  {\rm res}_g(a;F):=\frac{1}{(2\pi i)^k}\int_{\Gamma_a^\eps}
  \frac{F}{g_1\cdots g_k}dX_1\ldots dX_k,
\]
with $\Gamma_a^\eps=\{x\in U(a)\mid |g_i(x)|=\eps\}$, $U(a)$ a
neighbourhood of $a$ with $g^{-1}(0)\cap U(a)=\{a\}$ and $\eps$ so
small that $\Gamma_a^\eps$ lies relatively compact in $U(a)$.
$\Gamma_a^\eps$ is smooth for almost all $\eps$ by Sard's Theorem and has a
canonical orientation by the $k$-form $d(\arg g_1)\wedge\ldots\wedge
d(\arg g_k)|\Gamma_a^\eps$.  (This local residue of course makes sense
for holomorphic $g$ and $F\in\calo_a$, but the polynomial case, to
which the general case may easily be reduced, is sufficient for our
purposes). We define the total residue
\[
  {\rm Res}_g(F):=\sum_{a\in g^{-1}(0)}{\rm res}_g(a;F),
\]
which is also known as Grothendieck residue symbol
${F\choose g_1,\ldots,g_k}$ in the context of duality
theory in algebraic geometry \cite{hartshorne}. Let
$J=\det\left(\frac{\di g_i}{ \di X_j}\right)$ be the
Jacobian of $g$. Then for regular values $y$ of $g$
\[
  {\rm Res}_{g-y}(F)=\sum_{x\in g^{-1}(y)}\left(\frac{F}{J}\right)(x)=
  \tr\left(\frac{F}{J}\right)(y).
\]
Therefore $\tr(F/J)$ extends holomorphically (surprise!) to a neighbourhood
of $0$ (the extension will be denoted $\tr(F/J)$ as well) and
\[
  {\rm Res}_g(F)=\tr\left(\frac{F}{J}\right)(0).
\]
One abstract feature in our setting is that
we have weights $d_i$ associated to $X_i$ and that our relations
$f_i^{[\omega]}$ form a standard basis of the relation ideal with respect to
these weights, i.e.\
\[
  ({\rm In}\skp f_1^{[\omega]},\ldots,{\rm In}\skp f_k^{[\omega]})=
  (f_1,\ldots,f_k)={\rm In}\skp (f_1^{[\omega]},\ldots,f_k^{[\omega]}),
\]
where ``${\rm In}\skp $'' means taking initial forms. This is trivial in our
case since we started with the homogeneous generators $f_1,\ldots,f_k$
of the relation ideal in a presentation of $H^*(M,\cz)$. In this
situation one can describe the residue map algebraically as follows:
\begin{prop}\label{residues}\sloppy
  Let $R=\cz[X_1,\ldots,X_k]/(g_1,\ldots,g_k)$ be artinian such that $\{g_i\}$
  is a standard basis of the relation ideal
  with respect to weights $d_i$ of $X_i$. Put
  $N:=\sum_i\deg g_i-\sum_i d_i$, $R_{<N}:=\left\{F\in\cz[X_1,\ldots,X_k]\
  |\ \deg F<N\right\}/(g_1,\ldots,g_k)$
  and $J=\det\left(\frac{\di g_i}{\di X_j}\right)$. Then
  \[
    R=R_{<N}\oplus\cz\cdot J,
  \]
  and the total residue map ${\rm Res}_g:\cz[X_1,\ldots,X_k]
  \rightarrow\cz$ factorizes via the projection onto the second factor
  as follows
  \[
    \cz[X_1,\ldots,X_k]\stackrel{{\rm can}}{\longrightarrow
    \hspace{-4.1ex}\longrightarrow}R
    \stackrel{{\rm pr}_2}{\longrightarrow}\cz\cdot J\longrightarrow\cz,
  \]
  where the last map sends $J$ to $\dim_\scz R$.
\end{prop}
\pf
To check the normalization we observe ${\rm
Res}_g(J)=\tr_g\Big(1\Big)(0)=$ degree of $g$ over $0=\dim_\scz R$ (the
latter equality is generally true by flatness if the covering space is
Cohen-Macaulay \cite{fischer}). Next it is well known that the
residue vanishes on elements of $(g_1,\ldots,g_k)$.  The claim thus
reduces to $\ker({\rm Res}_g)/(g_1,\ldots,g_k)=R_{<N}$.

If $\deg F<N$ then $F dX_1\ldots dX_k/g_1\cdots g_k$ extends to a
rational differential form $\ph$ on weighted projective space
$V=\pr_{(1,d_1,\ldots,d_k)}$ with polar divisor $D_1+\ldots+D_k\in{\rm
Div}(V)$, $D_i$ the natural extension of the divisor $(g_i)$ to $V$.
The point of course is that the divisor at infinity $V\setminus\cz^k$
is not a polar divisor of $\ph$. We claim $|D_1|\cap
\ldots\cap|D_k|\subset\cz^k$. In fact, the restriction of the
homogenization of $g_i$ to $V\setminus\cz^k\simeq
\pr_{(d_1,\ldots,d_k)}$ is just ${\rm In}\skp g_i$ and $V({\rm In}\skp g_1,
\ldots,{\rm In}\skp g_k) =\{0\}\in\cz^k$.  The latter follows because
otherwise $\dim{\rm Spec}\skp \cz[X_1,\ldots,X_k]/ ({\rm
In}\skp g_1,\ldots,{\rm In}\skp g_k)>0$ by homogeneity. But from the standard
basis property we have $\dim_\scz\cz[X_1,\ldots,X_k]/({\rm In}\skp g_1,\ldots,
{\rm In}\skp g_k)=\dim_\scz\cz[X_1,\ldots,X_k]/(g_1,\ldots,g_k)<\infty$. --- We
may thus desingularize $V$ (at infinity) without violating the
discreteness of $|D_1|\cap\ldots\cap|D_k|$ (we use the same notations
for the pulled-back objects). Now the global residue theorem tells that
on the compact manifold $V$ the sum of the local residues of $\ph$ with
respect to $D_1,\ldots,D_k$ equals zero (the local residue ${\rm
res}_g(a;F)$ is coordinate free in so far that it depends only on
the associated rational differential form
$FdX_1\ldots dX_k/{g_1\ldots g_k}$ and the divisors $(g_1),
\ldots,(g_k)$). This proves ${\rm Res}_g(F)=0$ in case $\deg
F<N$.

The second case is $F$ homogeneous of degree $>N$. We show $F\in({\rm
In}\skp g_1, \ldots,{\rm In}\skp g_k)$. This is an easy generalization
of a theorem of Macaulay (cf.\ e.g.\ \cite{tsikh}) to the weighted
situation. Namely, for any $G$ (wlog.\ homogeneous), setting $P=F\cdot
G$, $Q={\rm In}\skp (g_1) \cdots{\rm In}\skp (g_k)$, we have
$\frac{P}{Q}dX_1\ldots dX_k =(\deg P-\deg Q+\sum_i
d_i)^{-1}d\sigma=(\deg P-N)^{-1}d\sigma$ with
\[
  \sigma=\frac{P}{Q}\sum_{j=1}^k(-1)^{j-1}d_j\cdot X_j\,
  dX_1\ldots\widehat{dX_j}\ldots dX_k,
\]
where $\widehat{\hspace{3ex}}$ means that this entry is to be left out.
This is a simple check using the weighted Euler formula $\displaystyle
\sum_j d_j X_j\frac{\di H}{\di X_j}=\deg(H)\cdot P$, $H$ weighted
homogeneous (same proof as usual).  Thus
\[
  {\rm Res}_{{\rm In}\skp g}(F\cdot G)=\sum_a\int_{\Gamma_a^\eps}d\sigma=0
\]
for all (homogeneous) $G\in\cz[X_1,\ldots,X_k]$. But this implies $F\in
({\rm In}\skp g_1,\ldots,{\rm In}\skp g_k)$ (``duality theorem'', cf.\
\cite{tsikh}
--- this is Poincar\'e duality in our case!). Modulo $(g_1,\ldots,g_k)$
this means that we may reduce $F$ to lower degree. So proceeding by
induction $N$ turns out to be ``top-degree'' in $R$ in that all elements
of $R$ can be represented by polynomials of degree $\le N$.

What remains to be checked is that for $F$ homogeneous of degree $N$,
${\rm Res}_g(F)\neq0$ or $F/(g_1,\ldots,g_k)\in R_{<N}$. But in the
first instance ${\rm Res}_{{\rm In}\skp (g)}(F\cdot
G)=0\ \forall\ G\in\cz[X_1,\ldots,X_k]$ as shown above, and again
$F\in({\rm In}\skp g_1,\ldots,{\rm In}\skp g_k)$, i.e.\ modulo
$(g_1,\ldots,g_k)$, $F$ may be represented by a polynomial of degree
$<N$.
\qed
\begin{rem}\rm
The decomposition $R=R_{<N}\oplus\langle J\rangle$ is $\em not$ canonical but
rather depends on the particular presentation of $R$. This may be seen either
in elementary terms from the trannsformation formula for residues or as
manifestation of the choice of an isomorphism ${\rm Ext}^k_V(\calo_Z,
\Omega_V^k)\simeq H^0(V,\calo_Z)$ in the duality morphism
\[
  {\rm Ext}^k_V(\calo_Z,\Omega_V^k)\times H^0(V,\calo_Z)
  \longrightarrow\cz
\]
induced by the global residue. For quantum cohomology rings and
weightings coming from cohomology, however, $R_{<N}=\oplus_{d<N}H^d(M,\cz)$,
$N=\dim_\cz M$ and $\langle Y\rangle=H^{2N}(M,\cz)$, so the decomposition
{\em has} an invariant meaning in this case.
\qed
\end{rem}

In the quantum cohomology ring the top-degree class $\cz\cdot J$ is thus
spanned by the class $[\Omega]$ of the volume form. Let $F_{[\Omega]}$
be a polynomial of degree $N$ representing $[\Omega]$ (modulo
$(f_i^{[\omega]})$ or modulo $(f_i)$, this will yield the same result),
and set $c= 1/{\rm Res}\skp _{f^{[\omega]}}(F_{[\Omega]})$. By the
interpretation of $\langle F\rangle_0$ as coefficient of $[\Omega]$ of
$F(\hat X_1,\ldots, \hat X_k)\in H^*_{[\omega]}(M)$ we conclude (for
$H^*(M,\cz)$ a complete intersection)
\[
  \langle F\rangle_0=c\cdot{\rm Res}_{f^{[\omega]}}(F)=
  c\cdot\tr_{f^{[\omega]}}\left(\frac{F}{J}\right)(0).
\]
To incorporate the higher genus case we prove
\begin{lemma}
  Notations as in the proposition and $F\in\cz[X_1,\ldots,X_k]$ let
  $B_F$ be the bilinear form $R\times R\rightarrow\cz$, $(\alpha,\beta)
  \mapsto{\rm Res}_g(F\cdot G_\alpha\cdot G_\beta)$ (with $G_\alpha$,
  $G_\beta\in\cz[X_1,\ldots,X_k]$ representing $\alpha$,
  $\beta$ --- this is well-defined), $\eta=B_1$ (``intersection form'').
  Then
  \[
    \tr_\eta B_F={\rm Res}_g(F\cdot J).
  \]
\end{lemma}
\pf
We give a basis-free, algebraic proof (some might find the brute force
method by adapting a basis to $\eta$ more enlightening).

By definition $\tr_\eta B_F$ is the trace of the endomorphism
$\mu_F:R\rightarrow R$ of multiplication by (the class of) $F$.
Putting $Z={\rm Spec}\skp R$ we have $R=\oplus_{z\in g^{-1}(0)}\calo_{Z,z}$.
For $z\in g^{-1}(0)$, $F-F(z)$ is nilpotent in $\calo_{Z,z}$
(for it has value $0$ in $z$), so $\mu_{F-F(z)}$ has trace $0$ and
\[
  \tr(\mu_F|\calo_{Z,z})=\tr(\mu_{F-F(z)}|\calo_{Z,z})
  +\tr(\mu_{F(z)}|\calo_{Z,z})=F(z)\cdot\dim_\scz\calo_{Z,z}.
\]
Furthermore, $\dim_\scz\calo_{Z,z}=\deg_z(g)$, the local mapping degree
of $g$ at $z$, so summing up we get
\[
  \tr_\eta B_F=\sum_{z\in g^{-1}(0)}\deg_z(g)\cdot F(z)=\tr_g(F),
\]
which is nothing but ${\rm Res}_g(F\cdot J)$ as claimed.
\qed
\vspace{3ex}

In view of the reduction formula to lower genus stated above we conclude
a rigorous version of the ``handle gluing formula'', previously
established in QFT by Witten \cite{witten0}.
\begin{prop}
  Let $H^*(M,\cz)=\cz[X_1,\ldots,X_k]/(f_1,\ldots,f_k)$ (the commutative
  complete intersection case) and $f_i^{[\omega]}$ the induced
  relations in $\qrM$ as in Theorem~\ref{presentation}. Put
  $J=\det\Big(\frac{\di f_i^{[\omega]}}{\di X_j}\Big)$. Then
  \[
    \langle F\rangle_g=\langle J\cdot F\rangle_{g-1}
  \]
  holds for all $F\in\cz[X_1,\ldots,X_k]$.
\qed
\end{prop}
In other words, multiplication by $J$ acts as ``attaching a handle''
to our Riemann surface. Our considerations so far yield the main
result of the present chapter:
\begin{theorem}\label{main4}
  Assumptions as in the preceding proposition then for all $F\in\cz[X_1,\ldots,
  X_k]$ the following holds
  \[
    \langle F\rangle_g=c\cdot{\rm Res}_{f^{[\omega]}}(J^g\cdot F)
    =c\cdot \lim_{\stackrel{\scriptstyle y\rightarrow0}{
    \stackrel{\scriptstyle y\mbox{\scriptsize\ regular}}{\scriptstyle
    \mbox{\scriptsize value of }f^{[\omega]} } } }
    \tr_{f^{[\omega]}}\left(J^{g-1}\cdot F\right)(y),
  \]
  with $c=1/{\rm Res}_{f^{[\omega]}}(F_{[\Omega]})$, $F_{[\Omega]}$ a
  polynomial representing the class $[\Omega]\in\qrM$ of the normalized
  volume form.
\qed
\end{theorem}

We emphasize that for $g>0$ or $0\in\cz^k$ a regular value of $f^{[\omega]}$
the right-hand side has the form $\sum_\nu a_\nu F(y_\nu)$ with
$g^{-1}(0)=\{y_\nu\}$ and constants $a_\nu$ {\em independent} of $F$!

Note also that $0$ is a regular value of $f^{[\omega]}$ iff $\dim_\scz
\calo_{Z,z}=1$ for all $z\in(f^{[\omega]})^{-1}(0)$, $Z={\rm Spec}\skp
\cz[X_1,\ldots,X_k]/(f_1^{[\omega]},\ldots,f_k^{[\omega]})$, which is
if and only if
\[
  \dim_\scz H^*(M,\cz)=\sharp\left(f^{[\omega]}\right)^{-1}(0)
\]
(``$\ge$'' always).
\vspace{3ex}

For the case of $\grass$ we have a basis of Schubert classes $\{\lambda_1,
\ldots,\lambda_k\}$ parametrized by sequences $n-k\ge\lambda_1\ge\ldots
\ge\lambda_k\ge0$. The latter are in \mbox{1--1} correspondence with subsets
$\{\lambda_k+1,\lambda_{k-1}+2,\ldots,\lambda_1+k\}$ of $\{1,\ldots,n\}$,
so $\dim_\scz H^*(\grass,\cz)={n\choose k}$.

To find $n \choose k$ distinct elements in
$\left(f^{[\omega]}\right)^{-1}(0)$ we use a description of the
cohomology ring coming from the study of the corresponding
Landau-Ginzburg model in physics \cite{vafa0}: There is a function $W$
(``Landau-Ginzburg potential'') s.th.\ in the notations of
Theorem~\ref{qrGrass}
\[
  Y_{n+1-i}(X_1,\ldots,X_k)=\frac{\di W}{\di X_i},
\]
and then $\displaystyle f_i^{[\omega]}=\frac{\di W^{[\omega]}}{\di
X_i}$ with $W^{[\omega]}:=W+(-1)^ke^{-[\omega](L)}\cdot X_1$. $W$ has a
simple description in terms of Chern roots, i.e.\ after
composition with $\Sigma:\cz^k\rightarrow\cz^k$,
$(\underline{\lambda})=(\lambda_1,\ldots,\lambda_k)
\mapsto\Sigma(\underline{\lambda})=
(-\sigma_1(\underline{\lambda}),\sigma_2(\underline{\lambda})
\ldots,(-1)^k\sigma_k(\underline{\lambda}))$,
$\sigma_i$ the elementary symmetric polynomials:
$W\circ\Sigma=-\frac{1}{n+1}\sum_i \lambda_i^{n+1}$. Thus
\[
  W^{[\omega]}\circ\Sigma=-\sum_{i=1}^k\left(\frac{\lambda_i^{n+1}}{n+1}
  -(-1)^k e^{-[\omega](L)}\cdot\lambda_i\right).
\]
This all is an easy formal consequence of the algebraic relations
between Chern and Segre classes on one side and Chern roots on the
other, cf.\ \cite{bertram.etal} for a mathematical account.  Now
$\displaystyle\frac{\di W^{[\omega]}}{\di X_i}\Big(\Sigma
(\underline\lambda)\Big)=0$, $i=1,\ldots,k$, if $\displaystyle\frac{\di
W^{[\omega]}\circ\Sigma}{\di\lambda_i}(\underline{\lambda})=0$ and if
$\Sigma$ is non-degenerate in $(\underline{\lambda})$.  The latter
obviously is equivalent to $\lambda_i\neq\lambda_j$ for $i\neq j$.
Solving the second equation means
\[
  \lambda_i^n=(-1)^k e^{-[\omega](L)},\ \ \ i=1,\ldots,k
\]
which has $n\cdot(n-1)\cdots(n-k+1)$ solutions with distinct
$\lambda_i$.  The fiber of $\Sigma$ over the regular value
$\Sigma(\underline\lambda)$ consists of the $k!$ permutations of
$\{\lambda_1,\ldots,\lambda_k\}$, so we get precisely $n\choose k$
distinct elements of $(f^{[\omega]})^{-1}(0)$, as wanted.
\begin{theorem}\label{VafIntForm}
  {\rm\bf (Formula of Vafa and Intriligator)}
  (Notations as in Theorem~\ref{qrGrass}.) For any $[\omega]\in
  H^{1,1}(\grass)$ the set $C^{[\omega]}$ of critical points of
  $W^{[\omega]}$ is finite and all of these are non-degenerate.
  Moreover, for any $F\in\cz[X_1,\ldots,X_k]$ the following formula
  for the genus $g$ GRW-invariants of $\grass$ holds:
  \[
    \langle F\rangle^{[\omega]}_g=
    (-1)^{n+k(k-1)/2}\sum_{\underline{x}\in C^{[\omega]}}
    \det\left(\frac{\di^2 W^{[\omega]}}{\di X_i\di X_j}\right)^{g-1}
    \hspace{-3ex}(\underline{x})\cdot F(\underline{x}).
  \]
\vspace{-5ex}

\qed
\end{theorem}
\pf
What is left is a check of the normalization. This amounts to calculate
the residue of $[\Omega]=(-1)^{n-k} X_k^{n-k}$ with respect to $f^{[\omega]}$.
We set $a=(-1)^k e^{[\omega](L)}$. Using standard properties of residues
one gets
\[
  {\rm Res}_{f^{[\omega]}}\left((-1)^{n-k} X_k^{n-k}\right)
  =\frac{(-1)^n}{k!}{\rm Res}_{{\underline\lambda}^n-\underline a}
  \left(\prod_{i=1}^k\lambda_i^{n-k}\prod_{i<j}(\lambda_i-\lambda_j)^2\right),
\]
where ${\underline\lambda}^n-\underline
a=(\lambda_1^n-a,\ldots,\lambda_k^n-a)$.
Note also that $\prod_{i<j}(\lambda_i-\lambda_j)^2$ is the squared Jacobian of
$\Sigma$ and that $k!$ is the degree of $\Sigma$. Since terms of degree less
than $k(n-1)$ modulo $(\lambda_i^n-a)$ have vanishing residue
(Proposition~\ref{residues}) only the term $\prod_{i=1}^k\lambda_i^{k-1}$ from
the expansion of $\prod_{i<j}(\lambda_i-\lambda_j)^2$ contributes. The
coefficient
of this term is $(-1)^{k(k-1)/2}\cdot k!$ as one sees by writing $\prod_{i<j}
(\lambda_i-\lambda_j)^2$ as determinant of a product of a Vandermonde-matrix
with its transposed. What is left is a multiple of the Jacobian $n^k\prod
\lambda_i^{n-1}$ of ${\underline\lambda}^n-\underline a$, the residue of
which is known to be the degree $n^k$ of ${\underline\lambda}^n-\underline a$.
Putting everything together we get the claimed normalization.
\qed
\vspace{3ex}

In concluding let us comment on the connection with the mathematical
formulation
of this formula given by Bertram, Daskalopoulos and Wentworth
\cite{bertram.etal}.
They showed that for any Riemann surface $\Sigma$ the moduli space
$\calm(d,\Sigma)$ (denoted $\calm$ in the sequel) of holomorphic
maps $f:\Sigma\rightarrow\grass$, $f_*[\Sigma]=d\cdot[L]$ has the expected
dimension provided the degree $d$ is sufficiently large. Moreover, the
compactification of $\calm$ as Grothendieck quot scheme $\overline
\calm_Q(d,\Sigma)$ (denoted $\overline\calm$ in the sequel), which parametrizes
sheaf quotients $\calo_\Sigma^n\rightarrow\calf\rightarrow0$ with fixed
Hilbert-polynomial, is generically reduced with irreducible reduction; the
universal quotient sheaf $\tilde\calf$ on $\Sigma\times\calm$
has a {\em locally free} kernel $\tilde\cale$ extending ${\rm ev}^*
S$ as subsheaf of $\calo^n_{\Sigma\times\overline\calm}$, $S$ the
tautological bundle on $\grass$. In particular it makes sense to talk about
the Chern classes $c_i(\tilde\cale)$. Let $\iota:\{t\}\times\overline\calm
\rightarrow\Sigma\times\overline\calm$ be the inclusion for some $t\in\Sigma$.
Then for $i_1,\ldots,i_k$ with $\sum_\nu \nu i_\nu=\dim\calm=k(n-k)(1-g)+dn$
they declared
\[
  \langle X_1^{i_1}\cdots X_k^{i_k}\rangle_\Sigma
  :=\left(\iota^*c_1(\tilde\cale)^{\wedge i_1}\wedge\ldots\wedge
  \iota^*c_k(\tilde\cale)^{\wedge i_k}\right)[\,\overline\calm\,].
\]
Presumably this invariant does not in general coincide with the
GRW-invariant: As we tried to explain in Remark~\ref{algmod} the
GRW-invariants are computed from the Chern classes of $\tilde{\rm
ev}^*(S)$, where $\tilde {\rm ev}:\Gamma\rightarrow\grass$ is an
extension of ${\rm ev}:\Sigma\times \calm\rightarrow\grass$ to some
blow-up $\pi:\Gamma\rightarrow\Sigma\times\overline\calm$. The trouble
is that $\pi^*\tilde\cale$, though a locally free subsheaf of
$\tilde{\rm ev}^*S$ of the same rank $k$, it need not coincide with the
latter. For instance this always happens if $\tilde\calf$ has torsion
since then $\calo_\Gamma^n/\pi^*\tilde \cale\simeq\pi^*\tilde\calf$ has
torsion as well. Additional contributions thus come precisely from the
first $k$ Chern classes of the torsion sheaf $\tilde{\rm
ev}^*\calo(S)/\pi^*\tilde\cale$ (one should desingularize $\Gamma$ to
make sense of these Chern classes). A general vanishing theorem for
these seems to be unlikely.


\end{document}